\title{Edge-transitive embeddings of complete graphs}

\author{Gareth A. Jones\\
School of Mathematics\\
University of Southampton\\
Southampton SO17  1BJ, UK\\
{\tt G.A.Jones@maths.soton.ac.uk}\\
}

\documentclass[12pt]{article}
\usepackage[latin1]{inputenc}
\usepackage{latexsym}
\usepackage{amsmath}
\usepackage{amssymb}
\usepackage{amsfonts}
\usepackage{color}
\usepackage{tikz}
\usepackage{enumerate}
\newtheorem{thm}{Theorem}[section]

\newcommand{\Z}{\mathbb{Z}}
\newcommand{\M}{\mathcal{M}}
\newcommand{\F}{\mathbb{F}}

\newcommand{\K}{\mathcal{K}}

\date{}
\begin{document} 

\newpage

\maketitle

\begin{abstract}
Building on earlier work of Biggs, James, Wilson and the author, and using the Graver--Watkins description of the 14 classes of edge-transitive maps, we complete the classification of the edge-transitive embeddings of complete graphs.
\end{abstract}

\medskip


 \medskip
 
 
%%%%%%%%%%%%%%%%%

\section{Introduction}\label{intro}

A map on a surface is {\sl regular\/} (sometimes called {\sl fully regular}) or {\sl edge-transitive\/} if its automorphism group acts transitively on vertex-edge-face flags or edges respectively. A map on an orientable surface is {\sl orientably regular\/} if its orientation-preserving automorphism group acts transitively on arcs. The regular and orientably regular embeddings of complete graphs are all known, as are their edge-transitive orientable embeddings. Here we complete the classification of their edge-transitive embeddings by determining their non-regular embeddings in non-orientable surfaces.

In 1971 Biggs~\cite{Big}, building on earlier work of Heffter~\cite{Hef}, proved:

\begin{thm}
The complete graph $K_n$ has an orientably  regular embedding if and only if $n$ is a prime power.  \hfill$\square$
\end{thm}

The maps Biggs constructed to prove that this condition is sufficient are Cayley maps $\M_n(c)$ for the additive groups of finite fields $\F_n$; in each case the generating set is the multiplicative group $\F_n^*=\F_n\setminus\{0\}$, taken in the cyclic order $1, c, c^2, \ldots, c^{n-2}$ where $c$ is a primitive element of $\F_n$ (that is, a generator for the cyclic group $\F_n^*$).  

\medskip

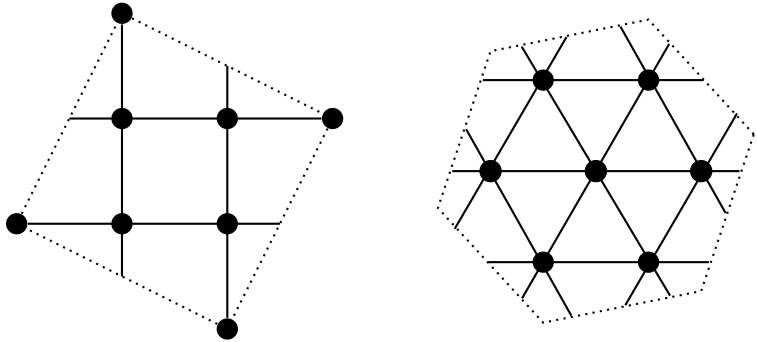
\begin{figure}[h!]

\begin{center}
 \begin{tikzpicture}[scale=0.7, inner sep=1mm]

\node (A) at (-6,3) [shape=circle, fill=black] {};
\node (B) at (-6,1) [shape=circle, fill=black] {};
\node (C) at (-4,1) [shape=circle, fill=black] {};
\node (D) at (-2,1) [shape=circle, fill=black] {};
\node (E) at (-8,-1) [shape=circle, fill=black] {};
\node (F) at (-6,-1) [shape=circle, fill=black] {};
\node (G) at (-4,-1) [shape=circle, fill=black] {};
\node (H) at (-4,-3) [shape=circle, fill=black] {};

\draw [thick] (-7,1) to (D);
\draw [thick] (E) to (-3,-1);
\draw [thick] (A) to (-6,-2);
\draw [thick] ((-4,2) to (H);

\draw[thick, dotted] (A) to (D);
\draw[thick, dotted] (D) to (H);
\draw[thick, dotted] (H) to (E);
\draw[thick, dotted] (E) to (A);

%%%%%%%

\node (a) at (3,0) [shape=circle, draw, fill=black] {};
\node (b) at (5,0) [shape=circle, draw, fill=black] {};
\node (c) at (2,1.73) [shape=circle, fill=black] {}; 
\node (d) at (4,1.73) [shape=circle, fill=black] {}; 
\node (e) at (2,-1.73) [shape=circle, fill=black] {}; 
\node (f) at (4,-1.73) [shape=circle, fill=black] {}; 
\node (g) at (1,0) [shape=circle, draw, fill=black] {};

\draw[thick, dotted] (6,0.7) to (4,2.88) to (1,2.28) to (0,-0.7) to (2,-2.88) to (5,-2.28) to (6,0.7); 

\node (h) at (5.88,0) {};
\node (H) at (0.12,0) {};
\draw [thick] (h) to (H);

\node (i) at (4.48,2.57) {};
\node (I) at (1.51,-2.57){};
\draw [thick] (i) to (I);

\node (j) at (1.5,2.52) {};
\node (J) at (4.5,-2.52) {};
\draw [thick] (j) to (J);

\node (k) at (0.7,1.73) {};
\node (K) at (5.2,1.73) {};
\draw [thick] (k) to (K);

\node (l) at (5.3,-1.73) {};
\node (L) at (0.8,-1.73) {};
\draw [thick] (l) to (L);

\node (m) at (3.38,2.9) {}; 
\node (M) at (5.55,-0.95) {}; 
\draw [thick] (m) to (M);

\node (n) at (2.62,-2.9) {}; 
\node (N) at (0.45,0.95) {}; 
\draw [thick] (n) to (N);

\node (o) at (5.77,1.25) {}; 
\node (O) at (3.47,-2.7) {}; 
\draw [thick] (o) to (O);

\node (p) at (0.23,-1.25) {}; 
\node (P) at (2.53,2.7) {}; 
\draw [thick] (p) to (P);

\end{tikzpicture}

\end{center}
\caption{The Biggs maps $\M_5(2)$ and $\M_7(3)$}
\label{K5K7}
\end{figure}

\noindent{\bf Example} For $n=4$ we have $\F_4=\{0, 1, c, c^2=c^{-1}=c+1\}$, and the corresponding map $\M_4(c)$ is the tetrahedral map on the sphere. The Biggs maps $\M_5(2)$ and $\M_7(3)$ are shown in Figure~\ref{K5K7}, with opposite sides of the outer square and hexagon identified to form the torus maps $\{4,4\}_{1,2}$ and $\{3,6\}_{1,2}$ in the notation of~\cite[\S8.3 and \S8.4]{CM}. Their mirror images are the maps $\M_5(3)\cong \{4,4\}_{2,1}$ and $\M_7(5)\cong \{3,6\}_{2,1}$.

\medskip

In 1985 James and the author~\cite{JJ} proved that the Biggs maps $\M_n(c)$ are the only orientably regular embeddings of complete graphs:

\begin{thm}\label{orKn}
A map $\M$ is an orientably regular embedding of $K_n$ if and only if $\M\cong\M_n(c)$ for some primitive element $c$ of $\F_n$. Moreover, $\M_n(c)$ and $\M_n(c')$ are isomorphic (as oriented maps) if and only if $c$ and $c'$ are equivalent under a field automorphism of $\F_n$.  \hfill$\square$
\end{thm}

It follows that if $n=p^e$ for some prime $p$ then there are, up to isomorphism, $\phi(n-1)/e$ orientably regular embeddings of $K_n$, one for each orbit of the Galois group of $\F_n$ (isomorphic to $C_e$, generated by the Frobenius automorphism $t\mapsto t^p$) on the $\phi(n-1)$ primitive elements of the field. The orientation-preserving automorphism group of $\M_n(c)$ is the affine group $AGL_1(\F_n)$; this map is regular if and only if $n=2, 3$ or $4$, in which case the full automorphism group is isomorphic to $V_4$, $D_6$ or $S_4$. 

If $n\ge 3$ the Petrie dual $P(\M_n(c))$ of $\M_n(c)$ is a non-orientable edge-transitive embedding of $K_n$, with the same automorphism group as $\M_n(c)$. 

In~\cite{Jam83} James classified the non-orientable regular embeddings of complete graphs (see also~\cite{Wil89} for an independent proof due to Wilson):

\begin{thm}\label{nonorKn}
The non-orientable regular embeddings of complete graphs $K_n$ are the maps $\{6, 2\}_3$, $\{4, 3\}_3$, $\{3, 5\}_5$ and $\{5, 5\}_3$ of characteristic $1, 1, 1$ and $-3$ for $n=3, 4, 6$ and $6$.  \hfill$\square$
\end{thm}

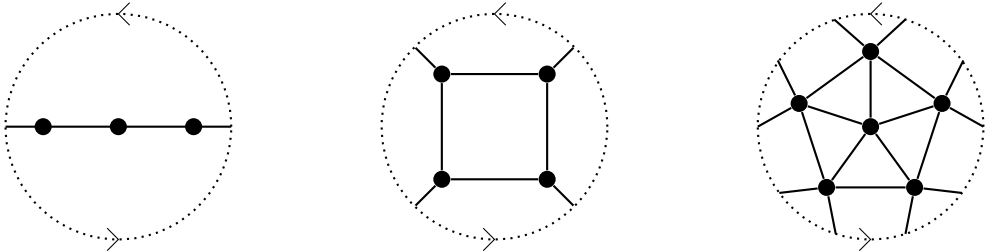
\begin{figure}[h!]
\begin{center}
\begin{tikzpicture}[scale=0.5, inner sep=0.8mm]

\draw [thick, dotted] (-10,0) arc (0:360:3);
\node (A) at (-15,0) [shape=circle, fill=black] {};
\node (B) at (-13,0) [shape=circle, fill=black] {};
\node (C) at (-11,0) [shape=circle, fill=black] {};
\draw [thick] (-16,0) to (-10,0);
\draw (-12.7,3.3) to (-13,3) to (-12.7,2.7);
\draw (-13.3,-3.3) to (-13,-3) to (-13.3,-2.7);

%%%%%%%

\draw [thick, dotted] (0,0) arc (0:360:3);
\node (a) at (-1.6,1.4) [shape=circle, fill=black] {};
\node (b) at (-4.4,1.4) [shape=circle, fill=black] {};
\node (c) at (-4.4,-1.4) [shape=circle, fill=black] {};
\node (d) at (-1.6,-1.4) [shape=circle, fill=black] {};
\draw [thick] (a) to (b) to (c) to (d) to (a);
\draw [thick] (a) to (-0.9,2.1);
\draw [thick] (b) to (-5.1,2.1);
\draw [thick] (c) to (-5.1,-2.1);
\draw [thick] (d) to (-0.9,-2.1);
\draw (-2.7,3.3) to (-3,3) to (-2.7,2.7);
\draw (-3.3,-3.3) to (-3,-3) to (-3.3,-2.7);

%%%%%%%

\draw [thick, dotted] (10,0) arc (0:360:3);
\node (a1) at (7,2) [shape=circle, fill=black] {};
\node (b1) at (8.9,0.618) [shape=circle, fill=black] {};
\node (c1) at (8.17,-1.616) [shape=circle, fill=black] {};
\node (d1) at (5.83,-1.616) [shape=circle, fill=black] {};
\node (e1) at (5.1,0.618) [shape=circle, fill=black] {};
\node (f1) at (7,0) [shape=circle, fill=black] {};
\draw [thick](a1) to (b1) to (c1) to (d1) to (e1) to (a1);
\draw [thick](a1) to (f1) to (b1);
\draw [thick](c1) to (f1) to (d1);
\draw [thick](e1) to (f1);

\draw [thick] (7.927,2.853) to (a1) to (6.033,2.853);
\draw [thick] (4.533,1.763) to (e1) to (4,0);
\draw [thick] (4.573,-1.763) to (d1) to (6.073,-2.853);
\draw [thick] (7.927,-2.853) to (c1) to (9.427,-1.763);
\draw [thick] (10,0) to (b1) to (9.467,1.763);

\draw (7.3,3.3) to (7,3) to (7.3,2.7);
\draw (6.7,-3.3) to (7,-3) to (6.7,-2.7);

\end{tikzpicture}

\end{center}
\caption{Regular embeddings of $K_n$ on the projective plane, $n=3, 4, 6$}
\label{projplKn}
\end{figure}

Here the notation $\{p,q\}_r$, from~\cite[\S8.6 and Table~8]{CM}, denotes the largest map of type $\{p,q\}$ with Petrie length $r$. The first three of these maps, on the real projective plane, are the antipodal quotients of a hexagon, a cube and an icosahedron on the sphere (see Figure~\ref{projplKn}, where antipodal boundary points of the discs are identified). The first two are the Petrie duals of the unique regular embeddings $\M_n(c)$ of $K_3$ and $K_4$ on the sphere, with the same automorphism groups $D_6$ and $S_4$, while the last two are a Petrie dual pair with automorphism group $PSL_2(5)\cong A_5$.

In~\cite{Jam90} James extended Theorem~\ref{orKn} to a classification of the orientable edge-transitive embeddings of $K_n$. If $3<n=p^e\equiv 3$ mod~$(4)$ where $p$ is prime, and $c$ is a primitive element of $\F_n$, let $\M_n(c,j)$ be the Cayley map for $\F_n$ with generating set $\F_n^*$, where now the cyclic ordering is
\[1, c^j, c^2, c^{j+2}, c^4, c^{j+4}, \ldots, c^{n-3}, c^{j+n-3}\]
for some odd element $j\in\Z_{n-1}\setminus\{1\}$. (Taking $j=1$ gives the orientably regular Biggs map $\M_n(c)$.) 

\begin{thm}\label{etransorKn}
A map $\M$ is an orientable edge-transitive embedding of $K_n$, which is not orientably regular, if and only if $\M\cong\M_n(c,j)$ for some $n$, $c$ and $j$ as above. As oriented maps, $\M_n(c, j)$ and $\M_n(c', j')$ are isomorphic if and only if $c$ and $c'$ are equivalent under ${\rm Gal}\,\F_n$ and $j'\equiv j$ or $2-j$ {\rm mod}~$(n-1)$.  \hfill$\square$
\end{thm}

These James maps $\M_n(c,j)$ have automorphism group $AHL_1(\F_n)$, the unique subgroup of index $2$ in $AGL_1(\F_n)$; their Petrie duals are non-orientable edge-transitive embeddings of $K_n$. The mirror image of $\M_n(c,j)$ is the  map $\M_n(c^{-1},2-j)\cong\M_n(c^{-1},j)$.

\medskip

\begin{figure}[h!]

\begin{center}
 \begin{tikzpicture}[scale=0.3, inner sep=1mm]
 
\draw [thick, dotted] (10,0) to (9.01,4.34) to (6.24, 7.82) to (2.23, 9.75) to (-2.23, 9.75) to (-6.24, 7.82) to (-9.01, 4.34) to (-10,0) to (-9.01,-4.34) to (-6.24,-7.82) to (-2.23, -9.75) to (2.23, -9.75) to (6.24, -7.82) to (9.01, -4.34) to (10,0);

\draw [thick] (1, -9.75) to (9.7, 1.1);
\draw [thick] (8.24, -5.30) to (5.19, 8.27);
\draw [thick] (9.29, 3.14) to (-3.23, 9.22);
\draw [thick] (-9.29, 3.14) to (3.23, 9.22);
\draw [thick] (-8.24, -5.30) to (-5.19, 8.27);
\draw [thick] (-1, -9.75) to (-9.7, 1.1);
\draw [thick] (7, -6.87) to (-6.9, -6.9);

\node (A) at (7.45,-1.8) [shape=circle, fill=black] {};
\node (B) at (6.04,4.7) [shape=circle, fill=black] {};
\node (C) at (0.106,7.67) [shape=circle, fill=black] {};
\node (D) at (-6.04,4.7) [shape=circle, fill=black] {};
\node (E) at (-7.45,-1.8) [shape=circle, fill=black] {};
\node (F) at (-3.24,-6.95) [shape=circle, fill=black] {};
\node (G) at (3.24,-6.95) [shape=circle, fill=black] {};

\draw [thick] (1,9.75) to (C) to (-1,9.75);
\draw [thick] (-7,6.87) to (D) to (-8.25,5.3);
\draw [thick] (-9.74,-1.19) to (E) to (-9.29,-3.14);
\draw [thick] (-5.15,-8.36) to (F) to (-3.34,-9.22);
\draw [thick] (5.15,-8.36) to (G) to (3.34,-9.22);
\draw [thick] (9.74,-1.19) to (A) to (9.29,-3.14);
\draw [thick] (7,6.87) to (B) to (8.25,5.3);

\node at (0,11) {$A$};
\node at (-8.9,-6.9) {$A$};
\node at (4.9,9.7) {$B$};
\node at (4.9,-9.7) {$B$};

\end{tikzpicture}

\end{center}
\caption{An edge-transitive embedding $\M_7(5,5)$ of $K_7$} \label{JamesK7}
\end{figure}
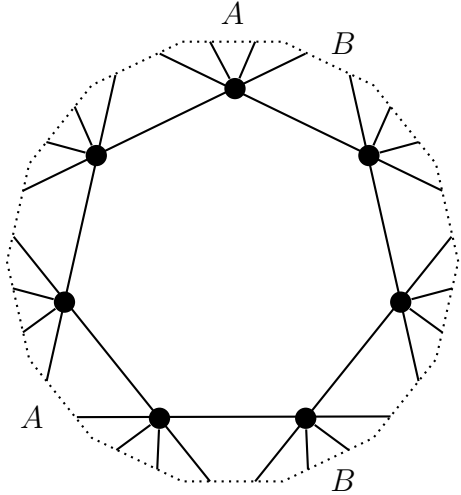

\noindent{\bf Example}  The James map $\M_7(5,5)$ is shown in Figure~\ref{JamesK7}, taken from~\cite{Jam90}. Two of the  identifications of sides of the outer $14$-gon are indicated by the letters $A$ and $B$; the others can be found by $C_7$ rotational symmetry. The vertices are identified with the elements $0, 1, \ldots, 6$ of $\F_7$, in clockwise order. There are seven triangular faces and three heptagons. The mirror image $\M_7(3,3)$ of this map is given by reflection in the horizontal axis. Each of these maps is a representation of the Fano plane, with the vertices as its points and the triples incident with its triangular faces forming its lines.

\begin{figure}[h!]

\begin{center}
 \begin{tikzpicture}[scale=0.3, inner sep=0.8mm]
 
\draw [thick, dotted] (10,0) to (9.01,4.34) to (6.24, 7.82) to (2.23, 9.75) to (-2.23, 9.75) to (-6.24, 7.82) to (-9.01, 4.34) to (-10,0) to (-9.01,-4.34) to (-6.24,-7.82) to (-2.23, -9.75) to (2.23, -9.75) to (6.24, -7.82) to (9.01, -4.34) to (10,0);

\draw [thick] (6.01,2.9) to (3.6,1.7);
\draw [thick] (1.48,6.5) to (0.89,3.9);
\draw [thick] (-4.16,5.21) to (-2.5,3.03);
\draw [thick] (-6.67,0) to (-4,0);
\draw [thick] (-4.16,-5.21) to (-2.5,-3.03);
\draw [thick] (1.48,-6.5) to (0.89,-3.9);
\draw [thick] (6.01,-2.9) to (3.6,-1.7);

\iffalse
\draw [thick] (1, -9.75) to (9.7, 1.1);
\draw [thick] (8.24, -5.30) to (5.19, 8.27);
\draw [thick] (9.29, 3.14) to (-3.23, 9.22);
\draw [thick] (-9.29, 3.14) to (3.23, 9.22);
\draw [thick] (-8.24, -5.30) to (-5.19, 8.27);
\draw [thick] (-1, -9.75) to (-9.7, 1.1);
\draw [thick] (7, -6.87) to (-6.9, -6.9);
\fi

\node (A) at (7.45,-1.8) [shape=circle, fill=black] {};
\node (B) at (6.04,4.7) [shape=circle, fill=black] {};
\node (C) at (0.106,7.67) [shape=circle, fill=black] {};
\node (D) at (-6.04,4.7) [shape=circle, fill=black] {};
\node (E) at (-7.45,-1.8) [shape=circle, fill=black] {};
\node (F) at (-3.24,-6.95) [shape=circle, fill=black] {};
\node (G) at (3.24,-6.95) [shape=circle, fill=black] {};

\node (a) at (7.45,1.8) [shape=circle, draw] {};
\node (b) at (6.04,-4.7) [shape=circle, draw] {};
\node (c) at (0.106,-7.67) [shape=circle, draw] {};
\node (d) at (-6.04,-4.7) [shape=circle, draw] {};
\node (e) at (-7.45,1.8) [shape=circle, draw] {};
\node (f) at (-3.24,6.95) [shape=circle, draw] {};
\node (g) at (3.24,6.95) [shape=circle, draw] {};

\iffalse
\draw [thick] (1,9.75) to (C) to (-1,9.75);
\draw [thick] (-7,6.87) to (D) to (-8.25,5.3);
\draw [thick] (-9.74,-1.19) to (E) to (-9.29,-3.14);
\draw [thick] (-5.15,-8.36) to (F) to (-3.34,-9.22);
\draw [thick] (5.15,-8.36) to (G) to (3.34,-9.22);
\draw [thick] (9.74,-1.19) to (A) to (9.29,-3.14);
\draw [thick] (7,6.87) to (B) to (8.25,5.3);
\fi

\node at (0,11) {$A$};
\node at (-8.9,-6.9) {$A$};
\node at (4.9,9.7) {$B$};
\node at (4.9,-9.7) {$B$};

\draw [thick] (B) to (6.01,2.9) to (a);
\draw [thick] (C) to (1.48,6.5) to (g);
\draw [thick] (D) to (-4.16,5.21) to (f);
\draw [thick] (E) to (-6.67,0) to (e);
\draw [thick] (F) to (-4.16,-5.21) to (d);
\draw [thick] (G) to (1.48,-6.5) to (c);
\draw [thick] (A) to (6.01,-2.9) to (b);

\draw [thick] (A) to (8,0) to (a);
\draw [thick] (B) to (5.1,6.2) to (g);
\draw [thick] (C) to (-1.8,7.8) to (f);
\draw [thick] (D) to (-7.3,3.6) to (e);
\draw [thick] (E) to (-7.3,-3.6) to (d);
\draw [thick] (F) to (-1.8,-7.8) to (c);
\draw [thick] (G) to (5.1,-6.2) to (b);

\draw [thick] (3.6,1.7) to (0.89,3.9) to (-2.5,3.03) to (-4,0) to (-2.5,-3.03) to (0.89,-3.9) to (3.6,-1.7) to (3.6,1.7);

\end{tikzpicture}

\end{center}
\caption{Vertices of $\M_7(5,5)$ and $\M_7(3,3)$ on Klein's surface}
\label{KleinK7}
\end{figure}
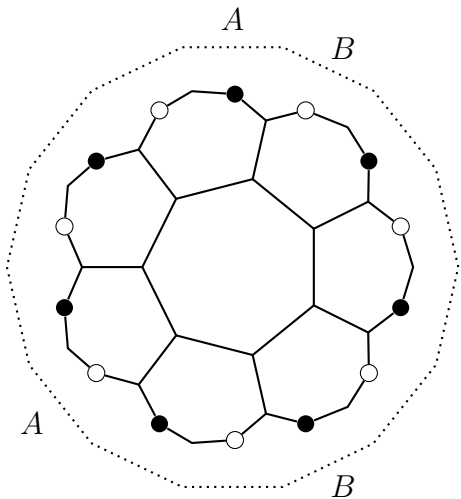

These two maps $\M$ can be drawn as follows on Klein's Riemann surface (or quartic curve) $\K$ of genus $3$ with automorphism group $PSL_2(7)$ (see~\cite{Lev}, for instance), so that ${\rm Aut}\,\M\;(\cong C_7\rtimes C_3)$ is contained in ${\rm Aut}\,\K$. Each of the eight Sylow $7$-subgroups $S\cong C_7$ of ${\rm Aut}\,\mathcal K$ preserves three of the $24$ faces of the ${\rm Aut}
\,{\mathcal K}$-invariant tessellation ${\mathcal T}=\{7,3\}_8$ of $\mathcal K$. These three faces are incident with 21 of the $56$ vertices of $\mathcal T$, and there are $21$ more vertices adjacent in $\mathcal T$ to these. This leaves $14$ vertices of $\mathcal T$, at graph-theoretic distance $2$ from the $S$-invariant faces, forming two orbits under $S$. The vertices in each orbit, joined pairwise by geodesics, determine a chiral pair of maps $\M\cong \M_7(5,5)$ and $\M_7(3,3)$, each invariant under the normaliser ${\rm Aut}\,\M=S\rtimes C_3$ of $S$ in ${\rm Aut}\,\mathcal K$. Figure~\ref{KleinK7} shows part of $\mathcal T$, with the central face invariant under the group $S$ generated by rotation through $2\pi/7$, and the orbits of $S$ on vertices of $\mathcal T$ yielding these two maps indicated in black and white. The vertices of each map are the centres of the triangular faces of the other. The three heptagonal faces of each map each contain one of the $S$-invariant faces of $\mathcal T$.

\medskip

In order to complete the classification of edge-transitive embeddings of complete graphs, it remains for us to deal with the non-regular non-orientable cases. The main result of this paper is as follows:

\begin{thm}\label{nonornonreg}
A map $\M$ is a non-orientable non-regular edge-transitive embedding of a complete graph $K_n$ if and only if $\M$ is isomorphic to the Petrie dual of a Biggs map $\M_n(c)$ for $n\ge 5$ or of a James map $\M_n(c,j)$ for $n\ge 7$.
\end{thm}

As an immediate corollary we have the following classification of the edge-transitive embeddings of complete graphs, showing that they are simply the maps listed above:

\begin{thm}\label{mainthm}
A  map $\M$ is an edge-transitive embedding of a complete graph $K_n$ if and only if $\M$ is isomorphic to a Biggs map $\M_n(c)$ or its Petrie dual, a James map $\M_n(c,j)$ or its Petrie dual, or one of the Petrie dual pair $\{3, 5\}_5$ and $\{5, 5\}_3$ of non-orientable regular embeddings of $K_6$. \hfill$\square$
\end{thm}

A more detailed description of these maps is given in Section~\ref{properties}.

%%%%%%%%%%%%%%%%%%%%%%%

\section{Algebraic theory of maps}\label{algthymaps}

Here we sketch the algebraic theory of maps developed in more detail elsewhere: see~\cite{JT}, for example, and~\cite{GT} for further background in topological graph theory. 

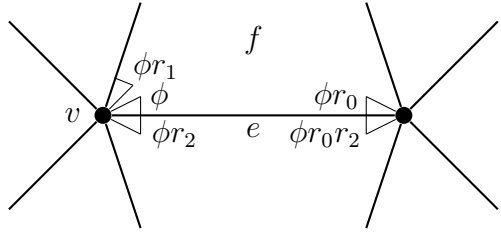
\begin{figure}[h!]
\begin{center}
\begin{tikzpicture}[scale=0.5, inner sep=0.8mm]

\node (c) at (0,0) [shape=circle, fill=black] {};
\node (d) at (8,0) [shape=circle, fill=black] {};
\draw [thick] (c) to (d);
\draw [thick] (c) to (1,-3);
\draw [thick] (c) to (1,3);
\draw [thick] (d) to (7,-3);
\draw [thick] (d) to (7,3);
\draw [thick] (c) to (-2.5,2.5);
\draw [thick] (c) to (-2.5,-2.5);
\draw [thick] (d) to (10.5,2.5);
\draw [thick] (d) to (10.5,-2.5);

\draw (c) to (1,0.5);
\draw (c) to (1,-0.5);
\draw (1,0.5) to (1,-0.5);
\draw (d) to (7,0.5);
\draw (d) to (7,-0.5);
\draw (7,0.5) to (7,-0.5);
\draw (c) to (0.8,0.8);
\draw (0.3,1) to (0.8,0.8);

\node at (-0.8,0) {$v$};
\node at (4,-0.4) {$e$};
\node at (4,2) {$f$};

\node at (1.5,0.5) {$\phi$};
\node at (6.2,0.5) {$\phi r_0$};
\node at (1.4,1.3) {$\phi r_1$};
\node at (1.9,-0.5) {$\phi r_2$};
\node at (5.9,-0.5) {$\phi r_0r_2$};

\end{tikzpicture}

\end{center}
\caption{Generators $r_i$ of $G$ acting on a flag $\phi=(v,e,f)$.} 
\bigskip
\label{flags}
\end{figure}

\begin{figure}[h!]
\begin{center}
\begin{tikzpicture}[scale=0.5, inner sep=0.8mm]

\draw [dashed] (-10,0) to (12,0);

\node (a) at (-8,3) [shape=circle, fill=black] {};
\draw [thick] (a) to (-8,0);
\draw (a) to (-7.5,2) to (-8,2);
\node at (-5.5,2.5) {$\phi r_0=\phi$};

\node (b) at (0,0) [shape=circle, fill=black] {};
\draw [thick] (b) to (0,3);
\draw (b) to (0.5,1) to (0,1);
\node at (2.5,1) {$\phi r_1=\phi$};

\node (c) at (7,0) [shape=circle, fill=black] {};
\draw [thick] (c) to (10,0);
\draw (c) to (8,0.5) to (8,0);
\node at (9,1.2) {$\phi r_2=\phi$};

\end{tikzpicture}

\end{center}
\caption{Flags fixed by $r_0, r_1$ and $r_2$.}
\label{fixedflags}
\end{figure}
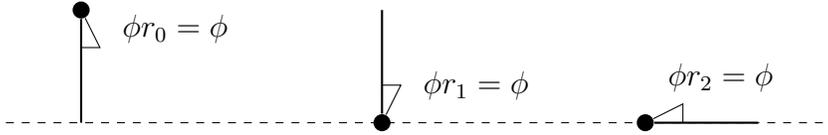

Each map $\mathcal M$ (possibly non-orientable or with non-empty boundary) determines a permutation representation of the group
\[\Gamma=\langle R_0, R_1, R_2\mid R_i^2=(R_0R_2)^2=1\rangle\cong V_4*C_2\]
on the set $\Phi$ of flags $\phi=(v,e,f)$ of $\mathcal M$, where $v, e$ and $f$ are a mutually incident vertex, edge and face. For each $\phi\in\Phi$ and each $i=0, 1, 2$, there is at most one flag $\phi'\ne \phi$ with the same $j$-dimensional components as $\phi$ for each $j\ne i$ (possibly none if $\phi$ is a boundary flag). Define $r_i$ to be the permutation of $\Phi$ transposing each $\phi$ with $\phi'$ if the latter exists, and fixing $\phi$ otherwise. (See Figures~\ref{flags} and~\ref{fixedflags} for the former and latter cases. In Figure~\ref{fixedflags} the broken line represents part of the boundary of the map.) Since $r_i^2=(r_0r_2)^2=1$ there is a permutation representation
\[\theta:\Gamma\to G:=\langle r_0, r_1, r_2\rangle\le{\rm Sym}\,\Phi\]
of $\Gamma$ on $\Phi$, given by $R_i\mapsto r_i$.

%%%%%%%%%%%%%

Conversely, any permutation representation of $\Gamma$ on a set $\Phi$ determines a map $\mathcal M$ in which the vertices, edges and faces are identified with the orbits on $\Phi$ of the subgroups $\langle R_1, R_2\rangle\cong D_{\infty}$,  $\langle R_0, R_2\rangle\cong V_4$ and  $\langle R_0, R_1\rangle\cong D_{\infty}$, incident when they have non-empty intersection.

The map $\mathcal M$ is connected if and only if $\Gamma$ acts transitively on $\Phi$, as we will always assume. In this case the stabilisers in $\Gamma$ of flags $\phi\in\Phi$ form a conjugacy class of subgroups $M\le\Gamma$, called {\em map subgroups}.

The map $\mathcal M$ is finite (has finitely many flags) if and only if $M$ has finite index in $\Gamma$, and it has non-empty boundary if and only if $M$ contains a conjugate of some $R_i$, or equivalently some $r_i$ has a fixed point in $\Phi$. In particular, $\M$ is orientable and with empty boundary if and only if $M$ is contained in the even subgroup $\Gamma^+$ of index $2$ in $\Gamma$, consisting of the words of even length in the generators $R_i$. We will assume unless stated otherwise that all maps considered have empty boundary.

The {\em automorphism group\/} $A={\rm Aut}\,{\mathcal M}$ of $\mathcal M$ is the centraliser of $G$ in ${\rm Sym}\,\Phi$. Then $A\cong N/M$ where $N:=N_{\Gamma}(M)$ is the normaliser  of $M$ in $\Gamma$. The map $\mathcal M$ is called {\em regular\/} if $A$ is transitive on $\Phi$, or equivalently $G$ is a regular permutation group, that is, $M$ is normal in $\Gamma$; in this case
$A\cong G\cong \Gamma/M$,
and one can identify $\Phi$ with $G$, so that $A$ and $G$ are the left and right regular representations of $G$. The map $\mathcal M$ is {\em edge-transitive\/} if $A$ acts transitively on its edges, or equivalently $\Gamma=NE$ where $E:=\langle R_0, R_2\rangle\cong V_4$.

The (classical) dual $D({\mathcal M})$ of $\mathcal M$ corresponds to the image of $M$ under the automorphism $\delta$ of $\Gamma$ fixing $R_1$ and transposing $R_0$ and $R_2$. The Petrie dual $P({\mathcal M})$ embeds the same graph as $\mathcal M$, but the faces are transposed with Petrie polygons, closed zig-zag paths turning alternately first right and first left at the vertices of $\mathcal M$; this operation corresponds to the automorphism $\pi$ of $\Gamma$ transposing $R_0$ with $R_0R_2$ and fixing $R_1$ and $R_2$. Both of these operations $D$ and $P$ preserve automorphism groups and regularity, but $D$ may change the embedded graph, and $P$ may change the underlying surface. The group $\Omega=\langle D, P\rangle$ of map operations generated by $D$ and $P$, introduced by Wilson in~\cite{Wil}, is isomorphic to $S_3$, permuting vertices, faces and Petrie polygons; it corresponds to the outer automorphism group ${\rm Out}\,\Gamma\cong {\rm Aut}\, E\cong S_3$ of $\Gamma$ acting on maps by permuting conjugacy classes of map subgroups~\cite{JT}.

%%%%%%%%%%%%%%%

\section{Edge-transitive maps}

\begin{figure}[h!]
\begin{center}
\begin{tikzpicture}[scale=0.25, inner sep=0.8mm]

\node (a) at (0,48) [shape=circle, fill=black] {};
\draw [very thick] (0,42) arc (-90:90:3);
\draw [dashed] (0,48) arc (90:270:3);
\node at (-5,45) {$1$};
% 1

%%%%%%%%%%%%%%%%%%%%%%%

\node (b) at (0,39) [shape=circle, fill=black] {};
\node (c) at (0,33) [shape=circle, fill=black] {};
\draw [very thick] (0,33) arc (-90:90:3);
\draw [dashed] (0,39) arc (90:270:3);
\node at (-5,36) {$2$};
% 2

%%%%%%%

\node (d) at (16,36) [shape=circle, fill=black] {};
\draw [dashed] (16,36) arc (0:360:3);
\draw [very thick] (d) to (10,36);
\node at (8,36) {$2^*$};
% 2*

%%%%%%%

\node (e) at (29,36) [shape=circle, fill=black] {};
\draw [dashed] (29,36) arc (0:360:3);
\draw [very thick] (e) to (26,36);
\node at (21,36) {$2^P$};
% 2P

%%%%%%%%%%%%%%%%%%%

\node (f) at (3,27) [shape=circle, fill=black] {};
\draw [very thick] (3,27) arc (0:360:3);
\node at (-5.5,27) {$2\,{\rm ex}$};
% 2ex

%%%%%

\node (g) at (13,27) [shape=circle, fill=black] {};
\draw [dashed] (16,27) arc (0:360:3);
\draw [very thick] (g) to (10,27);
\node at (7.5,27) {$2^*$ex};
% 2*ex

%%%%%%%

\node (h) at (27.5,27) [shape=circle, fill=black] {};
\draw [thin] (29,27) arc (0:360:3);
\draw [very thick] (h) to (24.5,27);
\node at (20.5,27) {$2^P$ex};
% 2Pex

%%%%%%%%%%%%%%%%%%%%

\node (i) at (-3,18) [shape=circle, fill=black] {};
\node (j) at (3,18) [shape=circle, fill=black] {};
\draw [very thick] (i) to (j);
\draw [dashed] (3,18) arc (0:360:3);
\node at (-5,18) {$3$};
% 3

%%%%%%%%%%%%%%%%%%%%

\node (k) at (0,9) [shape=circle, fill=black] {};
\node (l) at (3,9) [shape=circle, fill=black] {};
\draw [very thick] (k) to (l);
\draw [dashed] (3,9) arc (0:360:3);
\node at (-5,9) {$4$};
% 4

%%%%%%

\draw [dashed] (16,9) arc (0:360:3);
\node (m) at (10,9) [shape=circle, fill=black] {};
\draw [very thick] (m) to [out=30,in=90] (14,9);
\draw [very thick] (m) to [out=-30,in=-90] (14,9);
\node at (8,9) {$4^*$};
% 4*

%%%%%%

\draw [dashed] (23,12) to (29,12);
\draw [dashed] (23,6) to (29,6);
\draw [thick] [dotted] (23,12) to (23,6);
\draw [thick] [dotted] (29,12) to (29,6);
\node (n) at (23,6) [shape=circle, fill=black] {};
\node (o) at (29,12) [shape=circle, fill=black] {};
\draw [very thick] (n) to (o);
\draw [thin] (22.5,8.75) to (23,9.25);
\draw [thin] (23.5,8.75) to (23,9.25);
\draw [thin] (28.5,9.25) to (29,8.75);
\draw [thin] (29.5,9.25) to (29,8.75);
\node at (21,9) {$4^P$};
% 4P

%%%%%%%%%%%%%%%%%%%%%

\node (a) at (-1.5,0) [shape=circle, fill=black] {};
\node (b) at (1.5,0) [shape=circle, fill=black] {};
\draw [thin] (3,0) arc (0:360:3);
\draw [very thick] (a) to (b);
\node at (-5,0) {$5$};
% 5

%%%%%

\draw [thin] (16,0) arc (0:360:3);
\node (c) at (11,0) [shape=circle, fill=black] {};
\draw [very thick] (c) to [out=30,in=90] (15,0);
\draw [very thick] (c) to [out=-30,in=-90] (15,0);
\node at (8,0) {$5^*$};
% 5*

%%%%%%%

\node (d) at (23,0) [shape=circle, fill=black] {};
\node (e) at (29,0) [shape=circle, fill=black] {};
\draw [very thick] (29,0) arc (0:360:3);
\draw (25.75,3) to (26.25,3.5);
\draw (25.75,3) to (26.25,2.5);
\draw (26.25,-3) to (25.75,-3.5);
\draw (26.25,-3) to (25.75,-2.5);
\node at (21,0) {$5^P$};
% 5P

\end{tikzpicture}

\end{center}
\caption{The basic maps $\mathcal N(T)$ for the $14$ edge-transitive classes $T$}
\label{basicmaps}
\end{figure}
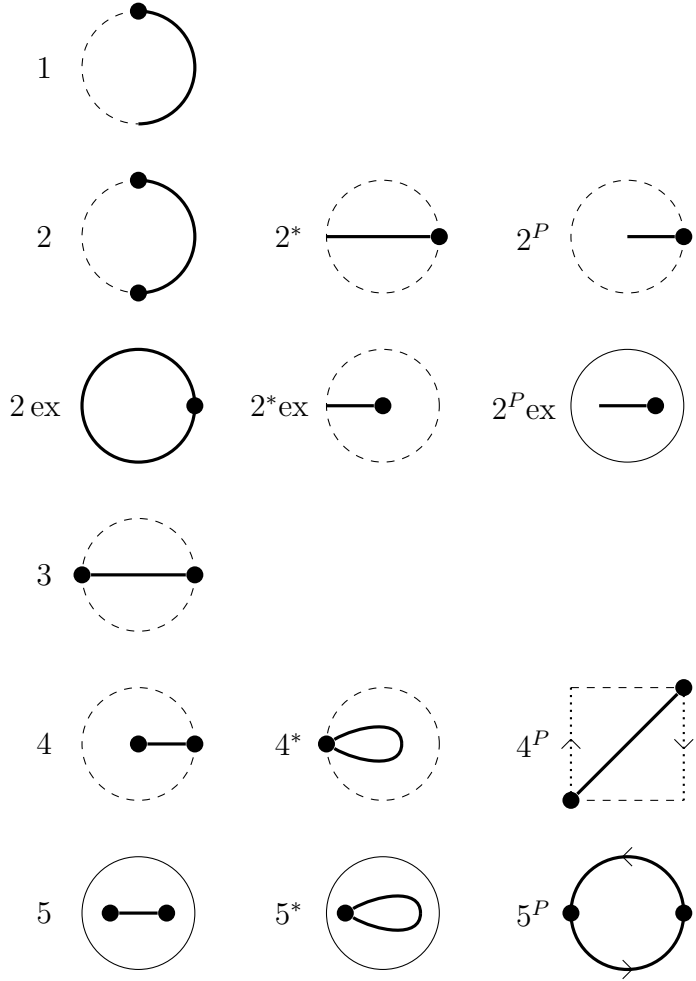

In 1997 Graver and Watkins~\cite{GW} partitioned the edge-transitive maps $\M$ into $14$ classes, distinguished by the isomorphism class of the quotient map ${\mathcal M}/{\rm Aut}\,{\mathcal M}$; in that year, Wilson~\cite{Wil97} gave a similar classification. These classes $T$ correspond to the $14$ isomorphism classes of maps ${\mathcal N}(T)$ with one edge, shown in Figure~\ref{basicmaps}, or equivalently to the $14$ conjugacy classes of {\sl  parent groups}, subgroups $N=N(T)$ of $\Gamma$ satisfying $\Gamma=NE$ (see~\cite[\S4]{Jon19}). The maps in Figure~\ref{basicmaps} are all on the closed disc, apart from ${\mathcal N}(2^P{\rm ex})$, ${\mathcal N}(5)$ and ${\mathcal N}(5^*)$ on the sphere, ${\mathcal N}(4^P)$ on the M\"obius band and ${\mathcal N}(5^P)$ on the real projective plane. The $14$ edge-transitive classes include class $1$, consisting of the regular maps, and class $2^P{\rm ex}$, consisting of the chiral (non-regular) orientably regular maps. Each map $\M$ in class $T$, with automorphism group $A$, is a regular covering by $A$ of the basic map ${\mathcal N}(T)$ for that class; it corresponds to a map subgroup $M\le\Gamma$ with $N_{\Gamma}(M)=N(T)$ and $N(T)/M\cong A$.

The six rows of Figure~\ref{basicmaps} correspond to the orbits of $\Omega$ on the $14$ classes: the duals of the maps in class $2$ form class $2^*$, while the Petrie duals of the latter form class $2^P$; similar remarks apply to the classes $2\,{\rm ex}, 4$ and $5$, whereas the classes $1$ and $3$ are invariant under $\Omega$.

The Reidemeister-Schreier process, applied to the inclusions $N(T)\le\Gamma$, gives the following presentations and free product decompositions:
\[N(1)=\Gamma=\langle R_0, R_1, R_2\mid R_i^2=(R_0R_2)^2=1\rangle\cong V_4*C_2,\]
\[N(2)=\langle S_1=R_1, S_2=R_1^{R_0}, S_3=R_2\mid S_1^2=S_2^2=S_3^2=1\rangle\cong C_2*C_2*C_2,\]
\[N(2\,{\rm ex})=\langle S_1=R_2, S=R_0R_1\mid S_1^2=1\rangle\cong C_2*C_{\infty},\]
\[N(3)=\langle S_0=R_1,\, S_1=R_1^{R_0},\, S_2=R_1^{R_2},\, S_3=R_1^{R_0R_2}\mid S_i^2=1\rangle\cong C_2*C_2*C_2*C_2,\]
\[N(4)=\langle S_1=R_1, S_2=R_1^{R_2}, S=(R_1R_2)^{R_0}\mid S_1^2=S_2^2=1\rangle\cong C_2*C_2*C_{\infty},\]
\[N(5)=\langle S=R_1R_2, S'=S^{R_0}\mid - \rangle\cong C_{\infty}*C_{\infty}\cong F_2.\]
(Here $F_2$ denotes a free group of rank $2$.) Applying elements of ${\rm Aut}\,E\cong\Omega$, permuting $R_0, R_2$ and $R_0R_2$, gives presentations for the other (isomorphic) parent groups in each orbit of $\Omega$. 

%%%%%%%%%%%%%%%%%%%

\section{Proof of Theorem~\ref{nonornonreg}}

It is easy to see that $P(\M_n(c))$ and $P(\M_n(c,j))$ have the stated properties.

Conversely, if $\M$ is a non-regular non-orientable edge-transitive embedding of $K_n$, then it is in one of the $10$ edge-transitive classes $T\ne 1$, $2^P{\rm ex}$, $5$ or $5^*$, since the maps in class~$1$ are regular, while those in the other three classes are orientable. The cases $n\le 3$ are easily dealt with, so assume that $n\ge 4$. If $T=2^*{\rm ex}$ or $5^P$ then $\M=P(\M')$ for some edge-transitive embedding $\M'$ of $K_n$ in class $2^P{\rm ex}$ or $5^*$; as such, $\M'$ is orientable, so it is isomorphic to a map $\M_n(c)$ or $\M_n(c,j)$ classified in Theorem~\ref{orKn} or \ref{etransorKn}, as required. We may therefore assume that $T\ne 2^*{\rm ex}$ or $5^P$. This leaves only $T=2$, $2^*$, $2^P$, $2\,{\rm ex}$, $3$, $4$, $4^*$ and $4^P$, so it is sufficient to eliminate these cases.

The edges of $K_n$ may be identified with the distinct pairs of vertices, so edge-transitivity implies that the group $A={\rm Aut}\,\M$ permutes these pairs transitively, that is, it acts $2$-homogenously on the vertices. It permutes the vertices faithfully since $n\ge 4$, and transitively since $K_n$ is not bipartite; the latter property rules out classes $T=2$, $3$ and $4$, since the corresponding maps ${\mathcal N}(T)$ have two vertices, so only the classes $T=2^*$, $2^P$, $2\,{\rm ex}$, $4^*$ and $4^P$ remain. If $A$ has odd order then, as a quotient of the group $N(T)$ with a free product decomposition given above, it must be cyclic; however, a cyclic group cannot be $2$-homogenous of degree $n\ge 4$, so $A$ has even order. It thus has an element transposing two vertices, so it is $2$-transitive on the vertices.

Each edge of $\M$ is therefore reversed by some automorphism (a reflection or half turn), so the single edge of ${\mathcal N}(T)$ must be free, ruling out the classes $T=2\,{\rm ex}$, $4^*$ and $4^P$. Only the classes  $2^*$ and $2^P$ remain, with $A$ sharply $2$-transitive on the vertices (since $\M$ is not regular), and the stabiliser of an edge generated by a reflection or a half-turn respectively, transposing its two incident vertices. Zassenhaus~\cite{Zas} showed that any sharply $2$-transitive finite group can be identified with $AGL_1(F)$ acting naturally on a near-field $F$. Since $\M$ is a map, the stabiliser $A_0$ of the vertex $0$ must be a cyclic or dihedral group of order $n-1$, acting regularly on the set $F\setminus\{0\}$ of neighbours of $0$. Now $A$ acts on the vertices as a Frobenius group $F\rtimes A_0$, so $A_0$ is a Frobenius complement; these contain at most one involution (see~\cite[Satz V.18.1(a)]{Hup} or~\cite[Theorem~18.1(iii)]{Pas}), so they cannot be dihedral, and hence $A_0$ is cyclic. Since $T=2^*$ or $2^P$, $N(T)$ is generated by involutions, and hence so are its epimorphic images $A$ and $A_0$, giving $n-1\le 2$, a contradiction. \hfill$\square$

%%%%%%%%%%%%%%%

\section{Summary of properties of the maps}\label{properties}

The following basic properties of the maps classified in Theorem~\ref{mainthm} are taken from~\cite{Big, Jam83, Jam90, JJ} or have been demonstrated earlier in this paper. 

For each prime power $n=p^e$ there are $\phi(n-1)/e$ Biggs maps $\M_n(c)$, all orientable. The mirror image of $\M_n(c)$ is $\M_n(c^{-1})$. For each $n=2, 3, 4$ the unique map $\M_n(c)$ is regular (in class~$1$), of type $\{2,1\}_2$, $\{3,2\}_6$ or $\{3,3\}_4$ and genus $0$, with automorphism group $V_4$, $D_6$ or $S_4$. For $n\ge 5$ these maps are orientably regular (in class $2^P{\rm ex}$) with automorphism group $AGL_1({\mathbb F}_n)$; they have type $\{m,n-1\}_{2p}$ where $m=(n-1)/2$ or $n-1$ as $n\equiv 3$ mod~$(4)$ or not, and have genus $(n^2-7n+4)/4$ or $(n-1)(n-4)/4$ respectively.

The maps $P(\M_n(c))$ for $n\ge 3$ are non-orientable, with the same automorphism group as $\M_n(c)$; they are in class $1$ or $2^*{\rm ex}$ as $n\le 4$ or $n\ge 5$. They have type $\{2p,n-1\}_m$ where $m$ is as above, and characteristic
\[\chi=n\left(1-\frac{n-1}{2}+\frac{n-1}{2p}\right).\]

If $3<n=p^e\equiv 3$ mod~$(4)$, there are $(n-3)\phi(n-1)/4e$ James maps $\M_n(c,j)$, all orientable, with mirror image $\M_n(c^{-1},2-j)$. They are in class $5^*$, with automorphism group $AHL_1(\F_n)$ consisting of the transformations $v\mapsto av+b$ in $AGL_1(\F_n)$ such that $a$ is a non-zero square in $\F_n$. This group has two orbits on the faces: if $j, 2-j\not\equiv (n-1)/2$ mod~$(n-1)$ they have cardinality $n(n-1,j)$ and $n(n-1,2-j)$ and the genus is
$\frac{n}{4}\left((n-3)-2(n-1,j)-2(n-1,2-j)\right)+1$, whereas if $j$ or $2-j\equiv (n-1)/2$ mod~$(n-1)$ they have cardinality $n$ and $n(n-1)/2p$ and the genus is
$(n-1)(n(p-1)-4p)/4p$.
There is a single orbit of $n(n-1)/l$ Petrie polygons, each polygon having length $l=2(n-1)/(n-1,2(j-1))$.

The maps $P(\M_n(c,j))$ are non-orientable, with automorphism group $AHL_1(\F_n)$; they are in class $5^P$, and have type $\{l,n-1\}$ and characteristic
\[\chi=n\left(1-\frac{n-1}{2}+\frac{n-1}{l}\right).\]

For $n=6$ the Petrie dual pair $\{3,5\}_5$ and $\{5,5\}_3$ are non-orientable, with $\chi=1$ and $-3$; they are regular, with automorphism group $PSL_2(5)\cong A_5$.

These maps $\M$ are all vertex-transitive. Indeed, in all cases ${\rm Aut}\,\M$ has a subgroup acting regularly on the vertices:  when $n$ is a prime power there is a unique such subgroup, namely the additive group of $\F_n$ acting by translations, whereas when $n=6$ there is a conjugacy classes of ten such subgroups, isomorphic to $S_3$. Apart from the James maps $\M_n(c,j)$, these maps are all face-transitive. They are all arc-transitive with the exception of the James maps and their Petrie duals.

%%%%%%%%%%%%%%%%%%%%

\section{Edge-transitive embeddings with boundary}\label{boundary}

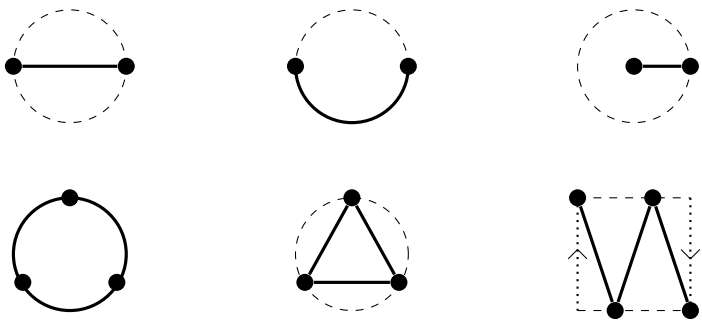
\begin{figure}[h!]
\begin{center}
\begin{tikzpicture}[scale=0.25, inner sep=0.8mm]

\node (a) at (-3,18) [shape=circle, fill=black] {};
\node (b) at (3,18) [shape=circle, fill=black] {};
\draw [very thick] (a) to (b);
\draw [dashed] (3,18) arc (0:360:3);

%%%%%%%

\node (c) at (12,18) [shape=circle, fill=black] {};
\node (d) at (18,18) [shape=circle, fill=black] {};
\draw [very thick]  (18,18) arc (0:-180:3);
\draw [dashed] (18,18) arc (0:180:3);

%%%%%%%

\node (e) at (33,18) [shape=circle, fill=black] {};
\node (f) at (30,18) [shape=circle, fill=black] {};
\draw [very thick] (e) to (f);
\draw [dashed] (33,18) arc (0:360:3);

%%%%%%%

\node (j) at (0,11) [shape=circle, fill=black] {};
\node (k) at (2.5,6.5) [shape=circle, fill=black] {};
\node (l) at (-2.5,6.5) [shape=circle, fill=black] {};
\draw [very thick] (3,8) arc (0:360:3);

%%%%%%%

\node (g) at (15,11) [shape=circle, fill=black] {};
\node (h) at (17.5,6.5) [shape=circle, fill=black] {};
\node (i) at (12.5,6.5) [shape=circle, fill=black] {};
\draw [very thick] (g) to (h) to (i) to (g);
\draw [dashed] (18,8) arc (0:360:3);

%%%%%%%

\draw [dashed] (27,11) to (33,11);
\draw [dashed] (27,5) to (33,5);
\draw [thick] [dotted] (27,11) to (27,5);
\draw [thick] [dotted] (33,11) to (33,5);
\node (m) at (27,11) [shape=circle, fill=black] {};
\node (n) at (29,5) [shape=circle, fill=black] {};
\node (o) at (31,11) [shape=circle, fill=black] {};
\node (p) at (33,5) [shape=circle, fill=black] {};
\draw [very thick] (m) to (n) to (o) to (p);
\draw [thin] (26.5,7.75) to (27,8.25);
\draw [thin] (27.5,7.75) to (27,8.25);
\draw [thin] (32.5,8.25) to (33,7.75);
\draw [thin] (33.5,8.25) to (33,7.75);

\end{tikzpicture}

\end{center}
\caption{Edge-transitive embeddings of $K_n$ with boundary, $n=2, 3$}
\label{Knbdy}
\end{figure}

As an addendum we show in Figure~\ref{Knbdy} those edge-transitive maps which embed complete graphs $K_n$ in surfaces with non-empty boundary. (A general theory of maps with boundary was developed by Bryant and Singerman in~\cite{BS}.)

There are three maps each for $n=2$ and $3$, and none for $n\ge 4$. (This follows easily from a wider study of edge-transitive maps with boundary in~\cite[\S17]{Jon19}.) The first five shown here are on the closed disc, while the last is on the M\"obius band. The first two embeddings of $K_2$ are in class~$1$, while the third is in class~$2$; their automorphism groups are isomorphic to $V_4$, $C_2$ and $C_2$. The three embeddings of $K_3$, all with automorphism group $S_3$, are in classes $1$, $2^*$ and $2^P$ respectively.

%%%%%%%%%%%%%%%%%%%%

\end{document}